\documentclass[a4article,12pt]{article}

\usepackage[T2A]{fontenc}
\usepackage[russian,english]{babel}
\usepackage[cp1251]{inputenc}
\usepackage{amsthm,amsmath}

\usepackage{amssymb}
\usepackage{graphicx}
\usepackage{color}

\makeatletter
\renewcommand{\@biblabel}[1]{#1.}
\makeatother

\theoremstyle{plain}
\newtheorem{theorem}{Теорема}

\theoremstyle{definition}

\newtheorem{corollary}{Следствие}

\newtheorem{example}{Пример}

\begin{document}
\selectlanguage{russian}
\begin{center}
\textbf{О ДЕФЕКТНЫХ ЧИСЛАХ ОПЕРАТОРОВ, ПОРОЖДЁННЫХ БЕСКОНЕЧНЫМИ ЯКОБИЕВЫМИ МАТРИЦАМИ}\\
И.Н.Бройтигам, К.А.Мирзоев\footnote{Работа первого автора выполнена при финансовой поддержке Министерства образования и науки РФ и Германской службы академических обменов (DAAD) по программе "Михаил Ломоносов" (ref: 325-A/13/74976)
и регионального  общественного фонда  содействия  отечественной
науке (грант  по  фундаментальным  проблемам  физики  и  математики), а второй автор поддержан  РНФ (грант № 14-11-00754).}
\end{center}

\noindent \textbf{1. Введение.}
Пусть $\mathbb{C}^m (m \ge 1)$ - евклидово $m$ - мерное пространство вектор-столбцов со стандартным скалярным произведением $y^*x= \sum\limits_{j=1}^m x_j \overline{y_j}$, где комплексные числа $x_j, y_j (j=1, \ldots, m)$ - координаты векторов $x$ и $y$ соответственно, а $*$ означает комплексно-сопряжённую матрицу. Пусть далее $A_j,B_j(j=0,1, \ldots)$ -  квадратные матрицы порядка $m$, элементы которых являются комплексными числами, причём  $A_j$ - самосопряжённые матрицы и  $B_j^{-1}$ существуют. Рассмотрим бесконечную якобиеву матрицу с матричными элементами
\begin{equation}
\label{JacobiM}
\mathbf{J}:=
\begin{pmatrix}
A_0 & B_0 & O & O & \ldots \\
B_0^* & A_1 & B_1 & O & \ldots \\
O & B_1^* & A_2 & B_2 & \ldots \\
\vdots & \vdots & \vdots & \vdots & \ddots
\end{pmatrix},
\end{equation}
где $O$ - нулевая матрица порядка $m$. Обозначим через $l_m^2$ гильбертово пространство бесконечных последовательностей $u=(u_0,u_1, \ldots), u_j \in \mathbb{C}^m$, со скалярным произведением $(u,v)=\sum\limits_{j=0}^{+\infty} v_j^* u_j.$  На всюду плотном в $l^2_m$ линейном многообразии финитных векторов отображение $u \rightarrow \mathbf{J}u$ определяет линейный симметрический, но незамкнутый оператор. Символами  $L$ и $D_L$ обозначим замыкание этого оператора и его область определения. Для $ u \in D_L$ оператор $L$ определяется формулой $Lu=lu$, где $lu=((lu)_0,(lu)_1, \ldots)$ и
$$
 (lu)_0:=A_0u_0 + B_0u_{1}, \;\; (lu)_j:=B^*_{j-1}u_{j-1} +A_ju_j +B_ju_{j+1}, \qquad j=1,2,\ldots,
$$
и является минимальным замкнутым симметрическим оператором, порождённым выражением $lu$ в пространстве $l^2_m$.

Подпространства  $\mathfrak{N}_-:=\{u \in l^2_m:L^*u=-iu \}
$ и $
\mathfrak{N}_+:=\{u \in l^2_m:L^*u=iu \}$
называются дефектными подпространствами (в точках $z=-i$ и  $z=i$), числа $n_-:=\mathrm{dim}\mathfrak{N}_-$ и $n_+:=\mathrm{dim}\mathfrak{N}_+$- дефектными числами, а пара $(n_-,n_+)$ - индексом дефекта оператора $L$. Дефектные числа удовлетворяют неравенствам $0 \le n_-,n_+ \le m$ и при этом если одно из них максимально, то и другое также максимально (см. \cite{Krein} - \cite{Kogan}). Известно, что для любой пары таких чисел $(n_-,n_+)$ существует  якобиевая  матрица $\mathbf{J}$ такая, что эта пара  является индексом дефекта
соответствующего симметрического оператора $L$ (см. \cite{Dukarev2}).

М.Г.Крейном (см. \cite{Krein},\cite{Krein1}) установлено, что всякой якобиевой матрице  $\mathbf{J}$ отвечает некоторая матричная степенная проблема моментов и эта проблема  имеет единственное (нормированное в некотором смысле) решение в том и только в том случае, когда одно из чисел  $n_-$ или $n_+$ равно нулю. Этот случай называется  определённым случаем матричной проблемы моментов.
Матричная проблема моментов называется вполне неопределённой, если $n_-=n_+=m$ (см.\cite{Krein},\cite{Krein1}). В этом случае иногда говорят, что  для матрицы $\mathbf{J}$ (для выражения $l$) имеет место вполне неопределённый случай.

В настоящей статье рассматриваются также бесконечные матрицы вида $J_m:=(c_{ij})$ с числовыми элементами $c_{ij} \in \mathbb{C},$ такими, что $c_{ij}=\overline{c_{ji}}$, $c_{ij}=0$, если $|i-j|>m$ и $c_{j,j+m} \ne 0$ $(i,j=0,1,\ldots)$. $J_m$ имеет, очевидно, тот же вид, что и $\mathbf{J}$ (см.(\ref{JacobiM})) с матричными элементами $A_j$ и $B_j$ размерности $m$, но теперь $B_j$ является нижней треугольной матрицей. В частности, при $m=1$ имеем $J_1=\mathbf{J}$, где $c_{jj}=A_j$ и $c_{j,j+1}=B_j$.
Если при этом   $A_j,B_j \in \mathbb{R}$ и $B_j > 0$, то  далее везде их будем обозначать  символами $a_j$ и $b_j$ $(j=0,1, \ldots)$ соответственно, а матрицу $J_1$ - символом $J$.
Отметим, что в этом случае речь идёт об определённости или неопределённости классической степенной проблемы моментов.

В данной работе исследуются вопросы максимальности, не максимальности и минимальности дефектных чисел оператора $L$ в терминах элементов якобиевой матрицы $\mathbf{J}$. Особое внимание уделено вопросу об условиях на числа $a_j$ и $b_j$ обеспечивающих реализацию случая определённости классической степенной проблемы моментов.

\noindent \textbf{2. Предварительные сведения.} Перечислим основные известные факты используемые в дальнейшем.  Справедлива следующая теорема (см. \cite{Kostuchenko_Mirzoev1}, \cite{Kostuchenko_Mirzoev3})
\begin{theorem}
\label{th1}
Для матрицы $\mathbf{J}$ имеет место вполне неопределённый случай тогда и только тогда, когда все решения векторного уравнения
\begin{equation}
\label{eq13}
(lu)_j=zu_j, \qquad j=1,2, \ldots,
\end{equation}
при $z=0$ принадлежат пространству $l_m^2$.
\end{theorem}

{\bf Замечание.} В формулировке этой теоремы в \cite{Kostuchenko_Mirzoev1} была допущена опечатка: нумерация условия (2) была указана, начиная с индекса $j=0$ вместо $j=1$. Но доказательство в \cite{Kostuchenko_Mirzoev1} проведено верно, в нём использовалась правильная нумерация и, кроме того, в формулировке этой теоремы в работе \cite{Kostuchenko_Mirzoev3} (теорема 0.1) опечатки уже нет. Воспринят опечатку за ошибку, автор работы \cite{Nagy} построил в $\S 7$ контрпример к теореме, которая получается из теоремы 1, если в условиях (2) нумерацию начинать с $j=0$. Исходя из этого, он указал на неверность результатов работ \cite{Kostuchenko_Mirzoev1} и \cite{Kostuchenko_Mirzoev3}, не заметив при этом, что в \cite{Kostuchenko_Mirzoev3} теорема 1 уже сформулирована верно. Таким образом, пример из $\S 7$ работы \cite{Nagy}, очевидно, не опровергает теоремы \ref{th1}.

Рассмотрим уравнение (\ref{eq13}) как матричное уравнение относительно неизвестной последовательности квадратных матриц $u_j$ размерности $m$ и обозначим через $P_j(z)$ и $Q_j(z)$ $(z \in \mathbb{C})$ его решения,
удовлетворяющие начальным условиям $P_0(z):=I$, $P_1(z):=B_0^{-1}(zI-A_0)$ и $Q_0(z):=O$, $Q_1(z):=B_0^{-1}$ соответственно. Функции $P_j(z)$ и $Q_j(z)$ являются многочленами $j$-го и $j-1$-го порядка от комплексного переменного $z$ с матричными коэффициентами и называются матричными многочленами первого и второго рода соответственно.

 Справедлива формула Кристоффеля-Дарбу (см. \cite{Berezanskiy}, Гл. VII, \S 2, п.11)
\begin{equation}
\label{KD}
(\overline{z}-z) \sum\limits_{j=0}^n P_j^*(z) P_j(z)=P_{n+1}^*(z)B_n^*P_n(z)-P_n^*(z)B_nP_{n+1}(z),
\end{equation}
В \cite{Kostuchenko_Mirzoev1} установлено, что двойная последовательность
$$
K_{ji}=Q_j(0)P^*_i(0)-P_j(0)Q_i^*(0), \qquad i,j=0,1,\ldots
$$
удовлетворяет матричному уравнению (\ref{eq13}) при $z=0$ по переменной $j$ при фиксированном $i$ и начальным условиям $K_{ii}=0$, $K_{i+1,i}=B_i^{-1}$ и справедлива
\begin{theorem}
\label{th2}
Для того чтобы для матрицы $\mathbf{J}$ имел место вполне неопределённый случай, необходимо и достаточно, чтобы для любой последовательности отрезков натуральных чисел $[n_k,m_k]$ таких, что $m_k \le n_{k+1} < m_{k+1}$ $(k=1,2, \ldots),$ выполнялось условие
$$
\sum\limits_{k=1}^{+\infty}\left(\sum\limits_{j=n_k}^{m_k} \sum\limits_{i=n_k}^j \|K_{ji}\|^2 \right )^{\frac 1 2} < +\infty,
$$
где $\| A \|$ означает норму матрицы $A$.
\end{theorem}
При $m=1$, очевидно, теорема \ref{th2} является критерием неопределённости классической степенной проблемы моментов.

Сформулируем одну теорему, относящуюся к случаю, когда оператор $L$ порождается якобиевой матрицей $J_m$.
\begin{theorem}
\label{th2_1}
Пусть
$$
\lim\limits_{j \to +\infty}\sup \frac{1}{\sqrt{c_{jj}^2+1}}\sum\limits_{k=1}^m    (|c_{j,j-k} |+|c_{j,j+k}|) <1.
$$
Тогда оператор $L$, порождённый матрицей $J_m$,  самосопряжён.
\end{theorem}
Эта теорема является частным случаем теоремы 2.2 работы \cite{Kostuchenko_Mirzoev3} (см. также \cite{Chistyakov}). \\
\textbf{3. Максимальность и не максимальность дефектных чисел оператора $L$.}
Используя теорему \ref{th2} можно доказать, что справедливы следующие две теоремы.
\begin{theorem}
\label{th2a}
Если
$$
\sum\limits_{j,i=1}^{\infty}  \|K_{ji} \| < \infty,
$$
то дефектные числа оператора $L$ максимальны.
\end{theorem}
\begin{theorem}
\label{th3}
Пусть при некотором $j$
\begin{equation}
\label{fth5}
\sum\limits_{n=1}^{\infty} \|K_{n+j,n} \| =\infty.
 \end{equation}
Тогда дефектные числа оператора $L$ не максимальны.
\end{theorem}
Заметим, что матрица $K_{ji}$ определяется из равенства
$$
K_{ji}=K_{ji}^0-\sum\limits_{k=i}^jK_{jk}^0A_kK_{ki}, \qquad i=0,1,\ldots, j=i+1,i+2, \ldots,
$$
где $$K_{i+2s,i}^0 =0 \,\,\,\mbox{и}\,\,\, K_{i+2s+1,i}^0=(-1)^sB_{i+2s}^{-1}B_{i+2s-1}^*B_{i+2s-2}^{-1}\ldots B^*_{i+1}B_i^{-1}, s,i =0,1, \ldots.$$
Учтя формулы для $K_{n+j,n}$ при $j=1,2, \ldots$ в (\ref{fth5}), получим иерархию достаточных признаков не максимальности дефектных чисел оператора $L$. В частности, при $j=1$ и $j=2$  имеем
$$K_{n+1,n}=B_n^{-1}, \qquad K_{n+2,n}=B_{n+1}^{-1}A_{n+1}B_n^{-1}$$
\noindent и поэтому справедливо
\begin{corollary}
\label{col1} \textit{Пусть  выполняется какое-либо из условий\\
\centerline{ a) $ \sum\limits_{n=1}^{\infty} \|B_n^{-1} \| =\infty,$ \qquad b) $ \sum\limits_{n=1}^{\infty} \|B_{n+1}^{-1}A_{n+1}B_n^{-1} \| =\infty.$} \\
Тогда дефектные числа оператора $L$ не максимальны.}
\end{corollary}
При $j=3$ и $j=4$, получим
$$K_{n+3,n}=-B_{n+2}^{-1}B_{n+1}^*B_n^{-1}+B_{n+2}^{-1}A_{n+2}B_{n+1}^{-1}A_{n+1}B_n^{-1},$$ $$K_{n+4,n}=B^{-1}_{n+3}B^*_{n+2}B^{-1}_{n+1}A_{n+1}B^{-1}_{n} + B^{-1}_{n+3}A_{n+3}B^{-1}_{n+2}B^{*}_{n+1}B^{-1}_{n}-B^{-1}_{n+3}A_{n+3}B^{-1}_{n+2}A_{n+2}B^{-1}_{n+1}A_{n+1}B^{-1}_{n}.$$
Используя далее теорему \ref{th3}, можно доказать, что справедливы следствия 2 и 3:
\begin{corollary}
\label{cor2} \textit{
Пусть\\
\centerline {$\sum\limits_{n=1}^{\infty} \| B_{n+2}^{-1}B_{n+1}^*B_n^{-1} \| < \infty \quad$
и
$\quad \sum\limits_{n=1}^{\infty} \|B_{n+2}^{-1}A_{n+2}B_{n+1}^{-1}A_{n+1}B_n^{-1} \|=\infty.
$}\\
Тогда дефектные числа оператора $L$ не максимальны. }
\end{corollary}
\begin{corollary} \textit{Пусть один из следующих рядов
 $$\sum\limits_{n=1}^{\infty} \|B^{-1}_{n+3}B^*_{n+2}B^{-1}_{n+1}A_{n+1} B^{-1}_{n} + B^{-1}_{n+3} A_{n+3}B^{-1}_{n+2}B^{*}_{n+1}B^{-1}_{n}\|,$$
  $$\sum\limits_{n=1}^{\infty} \| B^{-1}_{n+3}A_{n+3}B^{-1}_{n+2}A_{n+2}B^{-1}_{n+1}A_{n+1}B^{-1}_{n} \|$$
  сходится, а другой расходится. Тогда справедливо утверждение следствия \ref{cor2}.}
\end{corollary}
Теорема \ref{th2a} является, очевидно, признаком вполне неопределённости, а теорема \ref{th3} - признаком не вполне неопределённости якобиевой матрицы $\mathbf{J}$.

Отметим, что, в частности, при $m=1$ пункт $a)$ следствия \ref{col1} - это известный признак Карлемана определённости классической степенной проблемы моментов (\cite[гл. VII, \S1, теорема 1.3]{Berezanskiy}), а пункт $b)$ - признак Денниса и Уолла (\cite [гл. I, задача 2]{Ahiezer}).

\begin{example}
\label{pr1}
Пусть $1 \le p \le m-1$,  $J_p$ и $J_m$ (см. п.1) - якобиевые матрицы с элементами $c_{ij}^{\prime}$ и $c_{ij}^{\prime \prime}$ $(i,j=0,1,\ldots)$ соответственно, ненулевые элементы которых определяются равенствами $c_{j,j+p}^{\prime}=c_{j+p,j}^{\prime}=(j+1)^2$, $c_{j,j+m}^{\prime \prime}=c_{j+m,j}^{\prime \prime}=1$.  Применяя теорему 1 или \ref{th2a}, можно доказать, что индекс дефекта оператора $L$, порождённого матрицей $J_p$ равен $(p,p)$, а $J_m$ порождает, очевидно, самосопряжённый оператор в $l_m^2$. Поэтому индекс дефекта минимального оператора, порождённого суммой $J_p+J_m$ $(=:\mathbf{J}$) равен  $(p,p)$.
\end{example}
 Компоненты $A_j$ и $B_j^{-1}$  матрицы  $\mathbf{J}$ (см. (\ref{JacobiM})) из примера \ref{pr1}  являются матрицами размерности $m$ с целочисленными элементами, следовательно, они удовлетворяют и условию $a)$ и условию $b)$ следствия \ref{col1}. Таким образом, этот пример показывает, что при выполнении условий следствия \ref{col1} дефектные числа оператора $L$ могут принимать любые (равные) значения между нулём и $m-1$.

\noindent \textbf{4. Самосопряжённость оператора $L$.}
Известно несколько признаков самосопряжённости оператора $L$, порождённого якобиевой матрицей  $\mathbf{J}$ с матричными элементами (см., напр., \cite {Kostuchenko_Mirzoev3}, теорема 2.2). Как мы уже отмечали, в случае $m=1$ этот вопрос представляет особый интерес в связи с классической степенной проблемой моментов и при этом теорема \ref{th2} даёт, очевидно, необходимое и достаточное условие самосопряжённости оператора $L$, а теорема \ref{th3} демонстрирует эффективность этой теоремы. В недавней работе \cite{Velazquez} получен ряд новых признаков определённости классической степенной проблемы моментов. В этом пункте мы излагаем метод, позволяющий переносить часть результатов работы \cite{Velazquez} на случай оператора $L$, порождённого матрицей $\mathbf{J}$. Замечательно то, что здесь речь идёт не о немаксимальности дефектных чисел оператора $L$ как в теореме \ref{th3}, а о новых признаках его самосопряжённости.

Пусть $x \in \mathbb{C}^m$ и $x \ne 0$. Из формулы (\ref{KD})  можно извлечь (см. \cite[Гл. VII, \S 2, п.11, доказательство теоремы 2.9.] {Berezanskiy}), что справедливо неравенство
\begin{equation}
\label{BasKD}
\frac 1 {||B_n||} \le C \|P_{n+1}(i)x\| \|P_n(i)x\|,
\end{equation}
где $C={\| x \|^{-2}}$. С другой стороны, векторы $P_n(i)x$ удовлетворяют уравнению (\ref{eq13}) при $z=i$. Из этого и предположения о существовании матрицы $A_n^{-1}$ легко получить, что
\begin{equation}
\label{BasKD1}
\|P_n(i)x \| \le \| A_n^{-1}B_n \| \| P_{n+1}(i)x \| + \| A_n^{-1}B_{n-1}^* \| \| P_{n-1}(i)x \|,
\end{equation}
отсюда

$
\|P_n(i)x \| \le F_{1,n}(\| P_{n+1}(i)x \|+ \| P_{n-1}(i)x \|),
$
где $F_{1,n} = \| A_n^{-1}B_n \|  + \| A_n^{-1}B_{n-1}^* \|$.

Таким образом, применяя неравенство (\ref{BasKD}), получим
\begin{equation}
\label{Q1}
\frac 1 {||B_n||F_{1,n}} \le C \|P_{n+1}(i)x\|(\| P_{n+1}(i)x \|+ \| P_{n-1}(i)x \|).
\end{equation}
Из неравенства (\ref{BasKD1}) следует также, что $$\|P_n(i)x \| \le F_{2,n} (\|P_{n-2}(i)x\| + \|P_{n}(i)x\|  +\|P_{n+2}(i)x\|),$$
где $$F_{2,n}=\|A_n^{-1}B_{n-1}^*\|(\|A_{n-1}^{-1}B_{n-2}^*\|+\|A_{n-1}^{-1}B_{n-1}\|)+\|A_n^{-1} B_{n}\|(\|A_{n+1}^{-1}B_n^*\|+\|A_{n+1}^{-1}B_{n+1})\|).$$
Применяя снова (\ref{BasKD}), имеем
\begin{equation}
\label{Q2}
\frac 1 {||B_n||F_{2,n}} \le C \|P_{n+1}(i)x\| (\|P_{n-2}(i)x\| + \|P_{n}(i)x\|  +\|P_{n+2}(i)x\|).
\end{equation}
Продолжая аналогичные рассуждения и учитывая неравенства (\ref{BasKD}), (\ref{Q1}) и (\ref{Q2}), получим, что справедливо
\begin{equation}
\label{QQ}
\frac 1 {||B_n||F_{q,n}} \le C \|P_{n+1}(i)x\| \sum\limits_{k=0}^q \|P_{n-q+2k}(i)x\|, \quad q=0,1, \ldots,
\end{equation}
где $F_{0,n}=1$, $F_{1,n}$ и $F_{2,n}$ указаны выше, а при $q \ge 3$ выражения $F_{q,n}$ громоздки, поэтому мы их не приводим, отметим лишь, что они также как  $F_{1,n}$ и $F_{2,n}$ конструируются в соответствии с формулами, полученными в работе \cite{Velazquez} при $m=1$.

Исходя из неравенства (\ref{QQ}) можно доказать, что справедлива теорема
\begin{theorem}
\label{te1}
Оператор $L$, порождённый якобиевой матрицей $\mathbf{J}$ самосопряжён, если матричные элементы $A_n$ $(n=0,1, \ldots)$ обратимы и для какого-либо $q$ $(q=0,1, \ldots)$ выполняется условие
\begin{equation}
\label{QT}
\sum\limits_{n=q+1}^{\infty} \frac 1 {||B_n|| F_{q,n}} =\infty.
\end{equation}
 \end{theorem}

При $m=1$ эта теорема совпадает с теоремой 4.3 работы [11] и обобщает её на случай $m>1$.
Отметим, что в случае $m=1$  условие $a)$ следствия 1 и условие (\ref{QT}) при $q=0$ совпадают и являются признаком Карлемана самосопряжённости оператора $L$. Условие $b)$ следствия 1 и условие (\ref{QT}) при $q=1$ близки, но не совпадают даже при $m=1$.

\noindent \textbf{5. О степенях трёхдиагональных якобиевых матриц.} Пусть $k \in \mathbb{N}$ и $J$, как и выше,- якобиевая матрица с элементами $a_j$ и $b_j$ $(j=0,1,\ldots)$. Определим формальную степень порядка $k$ матрицы $J$, полагая как обычно $J^k=J(J^{k-1})$ и обозначим элементы этой матрицы символами $c^{(k)}_{ij}$ $(i,j=0,1,\ldots)$. Легко показать, что $c_{ij}^{(k)}=c_{ji}^{(k)}, c_{ij}^{(k)}=0$, если $|i-j|>k$ и $c_{j,j+k}^{(k)}$ отличен от нуля, т.е. согласно принятым во введении обозначениям, матрица $J^k$ является бесконечной матрицей вида $J_k$ и элементы $c_{n,n+s}^{(k)}$ при $k=1,2, \ldots, s=0,1, \ldots k$ и $n=0,1,\ldots$ определяются реккурентными соотношениями
$$
c_{n,n}^{(1)}=a_n, c_{n,n+1}^{(1)}=b_n,
$$
$$
c_{n,n+s}^{(k+1)}=b_{n-1}c_{n-1,n+s}^{(k)}+a_nc_{n,n+s}^{(k)}+b_nc_{n+1,n+s}^{(k)}, s=0,1, \ldots k-1,
$$
\begin{equation}
\label{Jk}
c_{n,n+k}^{(k+1)}=a_nc_{n,n+k}^{(k)}+b_nc_{n+1,n+k}^{(k)}, c_{n,n+k+1}^{(k+1)}=b_nc_{n+1,n+k+1}^{(k)}.
\end{equation}
Применяя теорему 1, можно доказать, что справедлива следующая
\begin{theorem}
Для якобиевой матрицы $J$ имеет место неопределённый случай тогда и только тогда, когда для матрицы $J^k$ имеет место вполне неопределённый случай.
\end{theorem}
Таким образом, для классической степенной проблемы моментов имеет место определённый случай тогда и только тогда, когда индексы дефекта минимального оператора $L(=:L_k)$, порождённого матрицей $J^k$ не максимальны.
В частности, если оператор $L_k$ самосопряжён, то оператор $L$, порождённый матрицей $J$ тоже является таковым. Используя это и применив теорему \ref{th2_1} можно доказать, что справедлива следующая теорема
\begin{theorem}
\label{th12}
Пусть $J$ - якобиевая матрица, $k \in \mathbb{N}$ и $c^{(k)}_{n,n+j}, (0 \le j \le k)$ - последовательность чисел, определяемая равенствами (\ref{Jk}). Тогда, если
$$
\lim\limits_{n \to \infty} \sup \frac 1 {\sqrt{\left|c_{nn}^{(k)}\right|^2+1}} \sum\limits_{j=1}^{k}  \left ( |{c_{n,n-j}^{(k)}}| + |{c_{n,n+j}^{(k)}}| \right )  <1,
$$
то оператор $L$, порождённый матрицей $J$ самосопряжён.
\end{theorem}
Полагая $k=2$, получаем следующее следствие из теоремы \ref{th12}.
\begin{corollary}
\label{col5}
Пусть $J$ - якобиевая матрица. Тогда если
\begin{equation}
\label{col_11}
\lim\limits_{n \to \infty} \sup  \frac {b_{n-1}b_{n-2} +|a_{n-1}+a_n|b_{n-1}+|a_n+a_{n+1}|b_n+b_nb_{n+1}} {b_{n-1}^2+a_n^2+b_n^2}<1,
\end{equation}
то оператор $L$, порождённый матрицей $J$ самосопряжён.
\end{corollary}
В заключении приведём один пример на применение следствия 4, который не охватывается теоремами 5 и 6.
\begin{example}
Пусть $a,b, \alpha \in \mathbb{R}$, $a \ne 0, b>0, \alpha >1$ и пусть $a_{n} \sim (-1)^n a n^{\alpha}$ и   $b_{n} \sim bn^{\alpha}$ при $n \to \infty$.  Тогда
предел в неравенстве (\ref{col_11}) равен $2b^2/(a^2+2b^2)$ и меньше 1. Следовательно, оператор $L$, порождённый матрицей $J$ самосопряжён.
Отметим, что если  $a_n=0$, то при условии $b_n \sim bn^{\alpha}$, $n \to \infty$ индексы дефекта оператора $L$ могут быть как $(0,0)$ так и $(1,1)$ (подробнее см. \cite[\S 1, п.3, пример 2]{Kostuchenko_Mirzoev3}).
\end{example}
Авторы благодарят профессора А.А. Шкаликова за проявленный интерес к работе и полезные замечания.

\noindent Бройтигам И.Н.,\\
 САФУ имени М.В. Ломоносова, \\
 Набережная Северной Двины, 17, \\
163002, г. Архангельск, Россия, \\
 e-mail:{irinadolgih@rambler.ru}\\
 Мирзоев К.А., \\
 МГУ имени М.В. Ломоносова, \\
 Ленинские Горы, 1, \\
 119991, г. Москва, Россия, \\
 e-mail:{mirzoev.karahan@mail.ru}
\end{document}